\documentstyle[amscd,amssymb]{amsart}

\newcommand{\seq}[3]{{#1}_{#2}, \ldots, {#1}_{#3}}
\newcommand{\Ct}{{\Bbb C}[[t]]}
\newcommand{\Cu}{{\Bbb C}[[u]]}
\newcommand{\Ctab}{{\Bbb C}[[t]][\alpha,\beta]}
\newcommand{\bl}{\frac{\partial}{\partial \l}}
\newcommand{\bs}{\frac{\partial}{\partial s}}
\newcommand{\bt}{\frac{\partial}{\partial t}}
\newcommand{\hS}{\hat{S}}
\newcommand{\hff}{HFF^*_g}
\newcommand{\ehff}{\widetilde{HFF}{}^*_g}
\newcommand{\rhff}{\overline{HFF}{}^*_g}
\newcommand{\Dws}{D^{(w,\Sigma)}}
\newcommand{\Y}{\Sigma \times {{\Bbb S}}^1}

\newcommand{\la}{\langle}
\newcommand{\ra}{\rangle}
\newcommand{\ima}{{\bf i}}

\newcommand{\surj}{\twoheadrightarrow}
\newcommand{\inc}{\hookrightarrow}
\newcommand{\ar}{\rightarrow}
\newcommand{\bd}{\partial}
\newcommand{\x}{\times}
\newcommand{\ox}{\otimes}
\newcommand{\iso}{\cong}

\newcommand{\point}{\text{pt}}
\newcommand{\CP}{{\Bbb C \Bbb P}}

\newcommand{\rk}{\text{rk}}

\newcommand{\Sym}{\text{Sym}}

\newcommand{\PD}{\text{P.D.}}

\newcommand{\odd}{\hbox{\scriptsize odd}}
\newcommand{\even}{\hbox{\scriptsize even}}

\newcommand{\cI}{{\cal I}}

\newcommand{\cV}{{\cal V}}

\renewcommand{\AA}{{\Bbb A}}
\newcommand{\CC}{{\Bbb C}}

\newcommand{\HH}{{\Bbb H}}

\newcommand{\QQ}{{\Bbb Q}}

\newcommand{\SS}{{\Bbb S}}
\newcommand{\ZZ}{{\Bbb Z}}

\renewcommand{\a}{\alpha}
\renewcommand{\b}{\beta}

\newcommand{\g}{\gamma}

\newcommand{\h}{\theta}
\renewcommand{\l}{\lambda}
\newcommand{\k}{\kappa}
\newcommand{\s}{\sigma}

\renewcommand{\o}{\omega}
\newcommand{\p}{\phi}
\newcommand{\q}{\psi}

\renewcommand{\S}{\Sigma}
\newcommand{\D}{\Delta}

\theoremstyle{plain}
\begingroup
\theorembodyfont{\sl}
\newtheorem{thm}{Theorem}
\newtheorem{cor}[thm]{Corollary}
\newtheorem{lem}[thm]{Lemma}
\newtheorem{prop}[thm]{Proposition}

\endgroup

\theoremstyle{definition}
\newtheorem{defn}[thm]{Definition}

\theoremstyle{remark}
\newtheorem{rem}[thm]{Remark}

\title{Donaldson invariants of non-simple type $4$-manifolds}

\author{Vicente Mu\~noz}
\date{September, 1999.}
\address{ Departamento de Matem\'aticas \\ 
Facultad de Ciencias \\ Universidad Aut\'onoma de Madrid \\ 
Ciudad Universitaria Cantoblanco \\ 28049 Madrid\\  Spain}
\email{vicente.munoz@@uam.es}

\begin{document}

\begin{abstract}
  We give the shape of the Donaldson invariants of a general 
  $4$-manifold with $b_1=0$ and $b^+>1$. The resulting expression 
  involves modular forms and matches the physical predictions. 
\end{abstract}

\maketitle

\section{Introduction}
\label{sec:intro}

Donaldson invariants for a (smooth, compact, oriented) $4$-manifold
$X$ with $b_1=0$, $b^+>1$ and odd (and with a homology orientation)
are defined as linear functionals~\cite{KM}
$$
  D^w_X: \AA(X)= \Sym^*(H_0(X) \oplus H_2(X)) \ar \CC,
$$
where $w \in H^2(X;\ZZ)$. As for the grading of
$\AA(X)$, the elements in $H_2(X)$ have degree $2$ and the point 
$x \in H_0(X)$ has degree $4$. 
Let $d_0=-w^2 -\frac32 (1+b^+)$, so that for homogeneous 
$z\in \AA(X)$, $D^w_X(z)$ is non-zero only if 
$\frac12 \deg z \equiv d_0 \pmod 4$.

Set $\wp=x^2-4 \in \AA(X)$ as in~\cite{adjuncti}. 
By definition, $X$ is of simple type 
if the condition $D^w_X(\wp z)=0$ is satisfied 
for any $z\in \AA(X)$. Also
$X$ is of finite type when there is some $n\geq 0$ such
that $D^w_X(\wp^n z)=0$, for any $z\in \AA(X)$. The order of
finite type is the minimum of such $n$. All $4$-manifolds 
with $b^+>1$ are of finite type~\cite{hff} and the order of finite
type is independent of $w\in H^2(X;\ZZ)$ 
(\cite[theorem 5]{basic}).
The question about the existence of non-simple type 
$4$-manifolds with $b_1=0$ and $b^+>1$ is still open.

In~\cite{FS}~\cite{KM} a structure theorem for the Donaldson invariants
of $4$-manifolds of simple type with $b_1=0$, $b^+>1$ is given. For
such $4$-manifold $X$, we have
\begin{equation}
  D_X^w(e^{tD+\l x})= \sum_K (-1)^{\frac{K\cdot w+w^2}{2}} 
  a_K e^{2\l+Q(tD)/2+K\cdot tD} + \ima^{-d_0} \sum_K
  (-1)^{\frac{K\cdot w+w^2}{2}} 
  a_K e^{-2\l-Q(tD)/2+\ima K\cdot tD}
\label{eqn:KM}
\end{equation}
where the sum runs over is the set of basic classes $K\in
H^2(X;\ZZ)$ and $a_K$ are rational numbers.

A conjecture on the shape of the Donaldson invariants of non-simple type 
$4$-manifolds is presented in~\cite{Kr}. In~\cite{Moore-Witten}
Moore and Witten give the shape of the invariants using 
physical arguments. This is expressed in terms of modular forms.
It is our intention to give a mathematically
rigourous proof of their formula.

Let $X$ be a general $4$-manifold with $b^+>1$ and $b_1=0$.
By~\cite[theorem 6]{basic}, we know that there is a (finite)
set of distinguished cohomology classes $K \in H^2(X;\ZZ)$
(called basic classes) such that
for any $w \in H^2(X;\ZZ)$ there are non-zero
polynomials $p_K \in \Sym^* H^2(X) \ox \QQ[\l]$ (dependent on $w$)
such that
\begin{equation}
   D^w_X(e^{tD+\l x})=e^{Q(tD)/2+2\l}\sum_K p_K(tD,\l) e^{K\cdot tD}+
   \ima^{-d_0} 
   e^{-Q(tD)/2-2\l}\sum_K p_K(\ima tD,-\l) e^{\ima K\cdot tD},
  \label{eqn:basic}
\end{equation}
for any $D\in H_2(X)$. 
The collection of classes $K$ is independent of
$w$ and $K$ are lifts to integral cohomology of ${w}_2(X)$.
In~\cite[definition 1.4]{adjuncti}, the order of finite type
of a particular basic class $K$ is defined as 
$d(K)=\max \{2n| \text{$K$ is a basic class for 
$D^w_X(\wp^n\bullet)$}\}$. Using
$$
  D^w_X(\wp^n e^{tD+\l x})=
  \left( \frac{\bd^2}{\bd \l^2} -4 \right)^n D^w_X(e^{tD+\l x}),
$$
we see that 
$d(K)= \max \{2n | \bl p_K \neq 0 \}=2 \deg p_K$,
where the degree of $p_K$ is considered 
with respect to the variable $\l$ (this will be done in
the sequel unless otherwise specified).
This number is independent of $w \in H^2(X;\ZZ)$.
Note that by~\cite[theorem 1.7]{adjuncti}, one has
$|K\cdot\S|+\S^2 +2d(K) \leq 2g-2$, for any embedded
$\S\inc X$ of genus $g$, if $\S^2>0$ or $\S^2=0$ and $\S$
is odd in homology.
Finally the order of finite type of $X$ is 
$1+\max\{d(K)|\text{$K$ basic class of $X$}\}$.
We aim to prove the following theorem on the structure of
the Donaldson invariants of $X$.

\begin{thm}
\label{thm:main}
  Define the following power series in $q=e^{2\pi\ima\tau}$
  in terms of classical modular and quasi-modular forms (see 
  section~\ref{sec:4} for a review)
\begin{align*}
  h(\tau) &=\frac12 \h_2(\tau)\h_3(\tau), \\
  V(\tau) &=\frac{-3e_1(\tau)}{h(\tau)^2}, \\
  T(\tau) &=-\frac{G_2(\tau)+e_1(\tau)/2}{h(\tau)^2}.
\end{align*}
  Then for any $4$-manifold $X$ with $b_1=0$, $b^+>1$ and odd, 
  there are power series $f_K(\tau)$ in $q$ 
  (uniquely determined up to the term $q^{d(K)/2}$), 
  for each basic
  class $K\in H^2(X;\ZZ)$, such that 
  for any $w\in H^2(X;\ZZ)$ we have 
 $$
  D^w_X(e^{D+\l x})= \left[q^{-d(K)/2} e^{V(\tau)\l +
  T(\tau)Q(D)+ K\cdot D/2h(\tau)} (-1)^{\frac{K\cdot w +w^2}{2}}
  f_K(\tau) \right]_{q^0} +
 $$
 $$
  +\ima^{-d_0} \left[q^{-d(K)/2} e^{-V(\tau)\l
  -T(\tau)Q(D)+ \ima K\cdot D/ 2h(\tau)} 
  (-1)^{\frac{K\cdot w +w^2}{2}}
  f_K(\tau) \right]_{q^0},
 $$
  where $[\cdot]_{q^0}$ stands for the coefficient of $q^0$ in
  the Laurent power series inside the bracket
  and $d_0=-w^2 -\frac32 (1+b^+)$.
  Moreover $f_{-K}(\tau)=(-1)^{(1+b^+)/2} f_K (\tau)$
  (at least modulo $q^{\frac{d(K)}{2}+1}$).
\end{thm}

\begin{rem}
\label{rem:1}
  If $X$ is of simple type then $d(K)=0$ for all basic
  classes and we recover the shape of the Donaldson invariants 
  given in~\eqref{eqn:KM}.
\end{rem}

\begin{rem}
\label{rem:12}
  A corollary to the proof of theorem~\ref{thm:main} 
  given in section~\ref{sec:5} is the
  following. Let $\tilde X=X\#\overline{\CP}^2$ be the blow-up 
  of a $4$-manifold $X$ with $b_1=0$, $b^+>1$ and odd, and
  let $E$ be the exceptional divisor. Then the basic classes
  of $\tilde X$ are exactly those classes $\tilde K=K\pm 
  (2n+1)E$
  with $K\in H^2(X;\ZZ)$ a basic class of $X$ and 
  $n\geq 0$, $n(n+1)\leq
  d(K)$. Moreover in this case $d(\tilde K)=d(K)-n(n+1)$.
  This agrees with the results in Seiberg-Witten 
  theory~\cite{FS-im}.
\end{rem}

\begin{rem}
\label{rem:physics}
  Theorem~\ref{thm:main} matches the physical predictions 
  in~\cite{Moore-Witten}. The following conjectures 
  remain to be proved: for every basic class $K$,
  $d(K)= {1\over 4}(K^2-2\chi-3\sigma)$ and
  $$
  f_K(\tau)=q^{-\frac{2\chi+3\s}{8}} \frac{SW(K)}{16} 
  \frac{\h_1(\tau)^{8+\s}}{a(\tau)h(\tau)}
  \left(2\frac{a(\tau)}{h(\tau)^2} \right)^{\frac{\chi+\s}{4}},
  $$
  where $a(\tau)=(4G_2(\tau)-e_1(\tau))/h(\tau)=16 q+\cdots$ 
  and $SW(K)$ is the Seiberg-Witten invariant 
  (see~\cite[equation (7.18)]{Moore-Witten}). One may
  check that $f_K(\tau)= 2^{1+\frac{7\chi+11\s}{4}}SW(K)+\cdots$,
  so that for simple type $4$-manifolds one recovers the
  usual conjecture.

  Our method of proof does not give more information on 
  $f_K(\tau)$. This is due to the fact that including any
  $z\in \CC[x]\subset \AA(X)$, one can run the same proof for
  $D^w_X(ze^{D+\l x})$ to get an expression analogous to
  that of theorem~\ref{thm:main}. On the other hand, once
  we have the result of theorem~\ref{thm:main},
  we always can find $z\in \CC[x]$ to arrange the 
  coefficients of $f_K(\tau)$ up to the term $q^{d(K)/2}$
  as we please.

  It has been pointed out to the author by G\"ottsche 
  that the results in~\cite{GZ} are the
  analogue for $4$-manifolds with $b^+=1$ 
  to our results for $4$-manifolds with $b^+>1$, but
  in~\cite{GZ} G\"ottsche and Don Zagier find out
  the formula for the analogue to $f_K$. The reason for
  this is that they use some homotopy invariance (that of
  the wall-crossing terms) that we do not have here at our
  disposal.
\end{rem}

We also rewrite theorem~\ref{thm:main} in the form suggested
in~\cite{Kr}. We have the following

\begin{thm}
\label{thm:main3}
  Let $K=K(k^2)$ and $E=E(k^2)$ be the elliptic integrals of
  first and second kind respectively, depending on the
  square of the modulus, $k^2$. Let $A(k^2)=\pi/2K(-\frac14
  k^2)$ and $B(k^2)=(2E(-\frac14 k^2)-K(-\frac14
  k^2))/K(-\frac14 k^2)$. Then
  for any $4$-manifold $X$ with $b_1=0$, $b^+>1$ and odd, 
  there are polynomials $P_K(\l)$ of degree $d(K)/2$, for
  every basic class $K\in H^2(X;\ZZ)$, such that for
  any $w\in H^2(X;\ZZ)$, we have 
  $$
  D^w_X(e^{D+\l x})= e^{2\l}\sum_K e^{B(\bl)Q(D)/2+
  A(\bl) K\cdot D} (-1)^{K\cdot w+w^2\over 2}
  P_K(\l) +\qquad
  $$
  $$ 
  \qquad +\ima^{-d_0} e^{-2\l}\sum_K 
  e^{-B(-\bl)Q(D)/2+A(-\bl)\ima K\cdot D} 
  (-1)^{K\cdot w+w^2\over 2} P_K(-\l),
  $$
  where $d_0=-w^2 -\frac32 (1+b^+)$. Moreover 
  $P_{-K}(\l)=(-1)^{(1+b^+)/2} P_K(\l)$.
\end{thm}

The paper is organized as follows. In section~\ref{sec:2} we 
study the Fukaya-Floer homology of the three-manifold
$Y=\Y$ and obtain some new relations which in
section~\ref{sec:3} are translated into a structure theorem 
for the Donaldson series $D^w_X(e^{tD+\l x})$ 
of a $4$-manifold $X$ with $b_1=0$ and $b^+>1$ and odd.
In section~\ref{sec:4} we review the modular forms
that we use in this paper and in section~\ref{sec:5}
we recast the expression obtained in section~\ref{sec:3}
into a $q^0$-coefficient of a power series in $q$. 
Using the universality of the expression thus obtained 
together with the general blow-up formula~\cite{FS-bl},
we get theorem~\ref{thm:main}. Finally in section~\ref{sec:6}
we derive theorem~\ref{thm:main3} from theorem~\ref{thm:main}.

{\em Acknowledgements:\/} The author is very grateful to 
Marcos Mari\~no
for teaching him on the interaction of modular forms and Donaldson
invariants and for explanations of the physical results.
He provided the author with the appropriate conjecture for
$f_K(\tau)$.
Also thanks to Lothar G\"ottsche for useful correspondence.

\section{Fukaya-Floer homology}
\label{sec:2}

Let us review the results of~\cite{hff} that we shall use in this
paper. For a (closed, connected, oriented) 
surface $\S=\S_g$ of genus $g \geq 1$,
consider the $3$-manifold $Y=\Y$. As in~\cite{hff}, $HFF^*_g=
HFF^*(Y,\SS^1)$ will stand for the Fukaya-Floer cohomology 
of $Y$ for the loop $\SS^1 \subset Y=\Y$ and $SO(3)$-bundle
with $w_2=\PD[\SS^1] \in H^2(Y;\ZZ/2\ZZ)$. 
In general we shall have the following situation: $X$ is a 
$4$-manifold with $b^+>1$, $b_1=0$ and containing an embedded
$\S \inc X$ of self-intersection zero and representing an
odd element in homology (so that there is $w\in H^2(X;\ZZ)$ with
$w\cdot \S \equiv 1\pmod 2$). Then we split $X=X_1\cup_Y A$, where
$A=A_g=\S \x D^2$ is a (closed) tubular neighbourhood of $\S$ and
$X_1$ is a $4$-manifold with boundary $\bd X_1=Y$. Let $D\in H_2(X)$
be a rational homology class which we represent as a $2$-cycle 
split as $D=D_1 +\D$, where $D_1\subset X_1$, $\bd D_1 =\SS^1\subset Y$
and $\D=\point\x D^2 \subset A$. Then for any $z \in \AA(X_1)=
\AA(\S^{\perp})\subset \AA(X)$ there is a relative 
Donaldson invariant~\cite{hff}
$$
  \p^w(X_1,z e^{tD_1}) \in HFF^*_g,
$$
such that 
\begin{equation}
  \Dws_X(z e^{tD})= \la\p^w(X_1,ze^{tD_1}),\p^w(A,e^{t\D})\ra,
\label{eqn:=}
\end{equation}
where $\Dws_X=D^w_X+D^{w+\S}_X$.
Let $\{\g_i\}$ be a symplectic basis for $H_1(\S;\ZZ)$ with 
$\g_i \cdot \g_{g+i}=1$, $1\leq i \leq g$. The elements
$\a= 2 \p^w(A,\S\,e^{t\D})$, $\b=-4\p^w(A,x\,e^{t\D})$ and
$\q_i= \p^w(A,\g_i\,e^{t\D})$ for $1\leq i \leq 2g$,
generate $\hff$ as a $\Ct$-algebra. %and $\g=-2 \sum \q_i\q_{g+i}$.
The product of $\hff$ has the property that
$\p^w(A, z_1 e^{t\D}) \p^w(X_1, z_2 e^{tD_1})= \p^w(X_1, z_1z_2 e^{tD_1})$,
for $z_1\in \AA(\S)$ and $z_2\in \AA(X_1)$.
For instance, $\a \p^w(X_1, z e^{tD_1})= \p^w(X_1,-2\S z e^{tD_1})$.

We ultimately want to extract information about the $4$-manifolds 
with $b^+>1$ 
and $b_1=0$, so it will be necessary to study the subspace of $\hff$ 
where the relative invariants $\p^w(X_1, z e^{tD_1})$ satisfying

  \hfill $(*)$ \parbox[t]{.8\linewidth}{$v=\p^w(X_1,z e^{tD_1})$ 
  where $X_1$ is a $4$-manifold with boundary $\bd X_1=Y=\Y$ such that  
  $X=X_1 \cup_Y A$ has $b_1=0$ and $b^+>1$, $z \in \AA(X_1)$, $D_1 
  \subset X_1$ is a $2$-cycle with $\bd D_1=\SS^1$ and $w|_Y=
  \PD[\SS^1]$} \hfill \break

\noindent live. We therefore give the following

\begin{defn}
\label{def:Vg} 
  Let $g\geq 1$. We define $\cV_g$ as the sub-$\Ct$-module of $\hff$
  generated by the relative invariants $v=\p^w(X_1,z e^{tD_1})$ 
  satisfying $(*)$.
\end{defn}

Clearly $\cV_g$ is a free $\Ct$-module, since $\hff$ is free as 
$\Ct$-module (see~\cite{hff}). Moreover $\cV_g$ is a $\Ctab$-algebra, 
since $\q_i\p^w(X_1,z e^{tD_1})=\p^w(X_1,\g_i z e^{tD_1})=0$, as
$\g_i=0$ in $H_1(X_1)=H_1(X)=0$. 

\begin{prop}
\label{prop:Vg}
  We have $\cV_g=\bigoplus\limits_{r=-(g-1)}^{g-1} R_{g,r}$ where
  $$
   R_{g,r}=\left\{ \begin{array}{ll}
   \Ctab/((\b+8)^d,\a-(4r + 2t) - f_{g,r}(\b+8)), & \text{$r$ odd} \\
   \Ctab/((\b-8)^d,\a-(4r\ima-2t)-f_{g,r}(\b-8)),\qquad & \text{$r$ even}
   \end{array}\right.
  $$
  for some $f_{g,r}(u) \in \Ct [u]$ polynomial of degree $d-1$ with 
  no independent term, where $d >0$ is the rank of $R_{g,r}$ as 
  $\Ct$-module.
\end{prop}

\begin{pf}
  We recall that the effective Fukaya-Floer homology $\ehff\subset \hff$
  studied in~\cite[section 5.4]{hff} is the sub-$\Ct$-module of $\hff$ 
  defined as $\cV_g$ but dropping the condition that 
  $X=X_1\cup_Y A$ has $b_1=0$. So $\cV_g\subset \ehff$.
  From~\cite[theorem 5.13]{hff} the eigenvalues of $(\a,\b)$ on 
  $\ehff$ are 
\begin{equation}
\label{eqn:eigenv}
  (4r+2t,-8),\text{ $r$ odd}, \quad (4r\ima-2t, 8), 
  \text{ $r$ even}, \qquad \text{with }-(g-1) \leq r \leq g-1.
\end{equation}
  As $\cV_g\subset \ehff$, the eigenvalues of $(\a,\b)$ on $\cV_g$
  are a subset of~\eqref{eqn:eigenv}, so we have the decomposition
  $\cV_g=\bigoplus\limits_{r=-(g-1)}^{g-1} R_{g,r}$, where
  in $R_{g,r}$,  $\a-(4r + 2t)$ and  $\b+8$ are nilpotent for $r$ odd,
  $\a-(4r\ima-2t)$ and $\b-8$ are nilpotent for $r$ even. 
  In principle it may happen that some $R_{g,r}$ is zero. 

  On the other hand, the reduced Fukaya-Floer homology
  $$
  \rhff= \ker(\b^2-64) \cap \bigcap_{i=1}^{2g} \ker \q_i \subset \hff
  $$
  is determined in~\cite[section 5.3]{hff} to be 
  $\rhff=\bigoplus\limits_{r=-(g-1)}^{g-1} {\bar R}_{g,r}$,
  where $\bar R_{g,r}$ are free $\Ct$-modules of rank $1$ such that 
  for $r$ odd, $\a=4r +2t$ and $\b=-8$ in $\bar R_{g,r}$. For $r$ even, 
  $\a=4r\ima -2t$ and $\b=8$ in $\bar R_{g,r}$. As
  $\ker(\b^2-64: \cV_g \ar \cV_g) \subset \rhff$ we have that
  $\ker(\b +(-1)^{r+1}8:R_{g,r}\ar R_{g,r}) \subset {\bar R}_{g,r}$. So
  $\ker(\b +(-1)^{r+1}8:R_{g,r}\ar R_{g,r})$ is either $\Ct$ or zero. 
  It is a simple algebra 
  exercise to check that this implies that $R_{g,r}$ is 
  as in the statement of the proposition, with $d\geq 0$ being the rank 
  of $R_{g,r}$.

  To check that in fact $R_{g,r} \neq 0$, for $-(g-1)\leq r \leq g-1$, 
  note that the vector $v\in \hff$ constructed 
  in~\cite[proposition 5.7]{hff} for $k=0$ actually lives in $\cV_g$, so
  $R_{g,\pm (g-1)}\neq 0$. Now the inequality~\eqref{eqn:0} below 
  (which holds without knowing a priori that all $R_{g,r}$,
  $-(g-1)\leq r \leq g-1$, are non-zero) gives the assertion.
\end{pf}

As a consequence of proposition~\ref{prop:Vg} the algebra $\cV_g$ is
cyclic, i.e.\ there is an epimorphism $\Ctab \surj \cV_g$. Therefore
$$
  \cV_g=\Ctab/\cI_g,
$$
for some ideal $\cI_g\subset \Ctab$. This ideal $\cI_g$ is the
ideal of relations of $\cV_g$, i.e.\ $f(\a,\b) \in \cI_g$ if and
only if $f(\a,\b) v=0$, for any $v=\p^w(X_1,ze^{tD_1}) \in \cV_g$.

\begin{rem}
\label{rem:determined}
  The Donaldson invariants of $4$-manifolds $X$ with $b_1=0$, $b^+>1$
  determine $\cV_g$. If $\p^w(X_1,z e^{tD_1})$ is in the conditions
  $(*)$ then we have the following chain of equivalences
$$
  \p^w(X_1,z e^{tD_1})=0 \iff 
  \la \p^w(X_1,z e^{tD_1}),\a^a\b^b \ra=0, \text{ for all $a, b\geq 0$}\iff 
$$
$$
  \iff \la \p^w(X_1,z e^{tD_1}),\p^w(A,\S^a x^b e^{t\D})\ra=0, 
         \text{ for all $a, b\geq 0$} \iff 
$$
$$
  \iff \la \p^w(X_1,z e^{tD_1}),\p^w(A,e^{t\D+s\S+\l x})\ra=0 \iff 
  \Dws_X(ze^{tD+s\S+\l x})=0,
$$
  where in the first equivalence we have used that $\a$, $\b$ and $\q_i$, 
  $1\leq i \leq 2g$ generate $\hff$ and $\q_i\p^w(X_1,z e^{tD_1})=0$, 
  for $1\leq i \leq 2g$.

  Thus $f(\a,\b) \in \cI_g$ if and only if
  $f(\a,\b) \p^w(X_1,ze^{tD_1})=\p^w(X_1,f(2\S,-4x)ze^{tD_1})=0$
  if and only if $\Dws_X(f(2\S,-4x)ze^{tD+s\S+\l x})=0$, for any 
  $4$-manifold $X$ with $b_1=0$ and $b^+>1$, an embedded
  $\S=\S_g \subset X$ with $\S^2=0$, $w\in H^2(X;\ZZ)$ with 
  $w\cdot \S\equiv 1\pmod 2$, $z\in\AA(\S^{\perp})$ and $D\in H_2(X)$
  with $D\cdot \S=1$.
\end{rem}

\begin{lem}
\label{lem:inc}
  $\cI_{g+1}\subset \cI_g$.
\end{lem}

\begin{pf}
  Let $f(\a,\b) \in \cI_{g+1}$. So $\Dws_X(f(2\S,-4x)ze^{tD+s\S+\l x})=0$
  for any $4$-manifold $X$ with $b_1=0$ and $b^+>1$, an embedded
  surface $\S=\S_{g+1} \subset X$ of genus $g+1$
  with $\S^2=0$, $w\in H^2(X;\ZZ)$ with 
  $w\cdot \S\equiv 1\pmod 2$, $z\in\AA(\S^{\perp})$ and $D\in H_2(X)$
  with $D\cdot \S=1$.

  Now if we have a $4$-manifold $X$ in the same situation with an embedded
  $\S=\S_g \subset X$ of genus $g$, one may add a trivial handle to $\S$ to
  obtain a new embedded surface of genus $g+1$ representing the same homology
  class. Therefore $\Dws_X(f(2\S,-4x)ze^{tD+s\S+\l x})=0$ and $f\in \cI_g$.
\end{pf}

  There is a natural epimorphism $\cV_{g+1} \surj \cV_g$.
  This yields in turn epimorphisms
  $R_{g+1,r}\surj R_{g,r}$, for any $-(g-1)\leq r \leq g-1$. Then
\begin{equation}
\label{eqn:0}
  \rk_{\Ct} R_{g+1,r}\geq \rk_{\Ct} R_{g,r}, 
\end{equation}
  for $g>|r|$. Let 
\begin{equation}
\label{eqn:00}
  x_r=\underset{g}{\sup} \> \rk_{\Ct} R_{g,r}. 
\end{equation}
  If $x_r<\infty$ then there is some $g_0$ such that $R_{g,r} \iso R_{g_0,r}$,
  for all $g\geq g_0$. If $x_r=\infty$ then the ranks of $R_{g,r}$ get
  bigger as $g\ar \infty$.

\begin{rem}
\label{rem:nilp}
  Let $v=\p^w(X_1,ze^{tD_1})\in \cV_g$ satisfying $(*)$. We may decompose
  $v=\sum v_r$, with $v_r\in R_{g,r}$. For any integer $\k <x_r$ we may
  find some $v$ such that $(\b+(-1)^{r+1}8)^{\k}v_r \neq 0$.
\end{rem}

\begin{cor}
\label{cor:Sr}
  There exist series $\hS_r(u,t) \in u\Ct [[u]]$ with 
  $$
   R_{g,r}=\left\{ \begin{array}{ll}
   \Ctab/((\b+8)^d,\a-(4r + 2t) - \hS_r(\b+8,t)), & \text{$r$ odd} \\
   \Ctab/((\b-8)^d,\a-(4r\ima-2t)-\hS_r(\b-8,t)), \qquad &\text{$r$ even}
   \end{array}\right.
  $$
  for any genus $g\geq 1$ and $-(g-1)\leq r\leq (g-1)$, where
  $d>0$ is the rank of $R_{g,r}$. The series $\hS_r$ is uniquely determined
  if $x_r=\infty$. If $x_r$ is finite then we impose that $\hS_r$ is
  a polynomial of degree less or equal than $x_r-1$. 
  With this extra condition $\hS_r$ is uniquely determined.
\end{cor}

\begin{pf}
  From  $R_{g+1,r}\surj R_{g,r}$ it is
  $\a-(4r+2t)-f_{g+1,r}(\b+8)=0$ in $R_{g,r}$. Therefore $f_{g+1,r}(u)$
  is equal to $f_{g,r}(u)$ plus (possibly) 
  terms of degrees strictly bigger than $\deg f_{g,r}$. 
  So there is a single series $\hS_r(u)$ agreeing with $f_{g,r}(u)$ 
  up to terms of degrees $\deg f_{g,r}$, for all $g>|r|$.
\end{pf}

\begin{rem}
\label{rem:Rr}
  By remark~\ref{rem:determined} we have that the series 
  $\hS_r$ are determined by the Donaldson invariants of
  $4$-manifolds with $b_1=0$ and $b^+>1$.
  If we write $R_r=\underset{g}{\varprojlim} R_{g,r}$ then 
  $x_r=\rk_{\Ct} R_r$. We have (e.g.\ for $r$ odd)
  $$
  R_r= \left\{ \begin{array}{ll} \Ctab/(\a-(4r + 2t) - \hS_r(\b+8)), & 
  \text{if $x_r=\infty$} \\ 
  \Ctab/((\b+8)^{x_r},\a-(4r + 2t) - \hS_r(\b+8)), \qquad  & 
  \text{if $x_r<\infty$} 
  \end{array} \right. 
  $$
\end{rem}

\section{Universal shape of the Donaldson invariants}
\label{sec:3}

Now we shall translate the knowledge on the Fukaya-Floer homology 
gathered in section~\ref{sec:2} to information on the Donaldson 
invariants of $4$-manifolds with $b_1=0$ and $b^+>1$. 
This section is devoted to proving the following result

\begin{thm}
\label{thm:main2}
  There are universal series $A(u), B(u) \in \Cu$ 
  of the form $1+\ldots$, such that
  for any $4$-manifold $X$ with $b_1=0$, $b^+>1$ and odd, 
  and any $w\in H^2(X;\ZZ)$, we have 
  $$
  D^w_X(e^{tD+\l x})= e^{2\l}\sum_K e^{B(\bl)Q(tD)/2+
  A(\bl) K\cdot tD} P_K(\l) +\ima^{-d_0} e^{-2\l}\sum_K 
  e^{-B(-\bl)Q(tD)/2+A(-\bl)\ima K\cdot tD} P_K(-\l),
  $$
  where $d_0=-w^2 -\frac32 (1+b^+)$. 
  The sum runs over the set of basic classes of $X$ and 
  $P_K(\l)$ are polynomials (dependent on $w$) of degree $d(K)/2$.
  Moreover $P_{-K}(\l)=(-1)^{d_0} P_K(\l)$.
\end{thm}

To start with, let us write
\begin{equation}
  S_r(u,t)=\left\{ \begin{array}{ll} \frac12 \hS_r(-4u,t), & \text{$r$ 
  odd} \\ \frac{\ima}2 \hS_{-r}(4u,\ima t), \qquad & \text{$r$ even}
\end{array} \right.
\label{eqn:Sr}
\end{equation}

\begin{lem}
\label{lem:inv}
  Suppose we are in the following situation

  \hfill $(**)$ \parbox[t]{.8\linewidth}{$X$ is a $4$-manifold with $b_1=0$ 
  and $b^+>1$ such that there
  is an embedded $\S\subset X$ with $\S^2=0$, $w\in H^2(X;\ZZ)$ with
  $w\cdot \S\equiv 1\pmod 2$, $z\in \AA(\S^{\perp})$ and $D\in H_2(X)$.}\hfill \break

  \noindent Then there exist polynomials $p_r(\l,t) \in 
\Ct[\l]$ such that  
$$
  D_X^w(ze^{tD+s\S+\l x})= e^{2\l}\sum_{r=-(g-1)}^{g-1}  
   e^{2rs +ts(D\cdot\S)} e^{S_r(\bl,t(D\cdot\S))s}p_r(\l,t) + 
$$
$$ 
   +\ima^{-d_0} e^{-2\l}\sum_{r=-(g-1)}^{g-1} e^{2r\ima s-ts(D\cdot\S)} 
   e^{S_r(-\bl,\ima t(D\cdot\S))\ima s} p_r(-\l,\ima t).
$$
  Note that $e^{S_r(\bl,t(D\cdot\S))s}p_r(\l,t)$ is polynomic in $s$ and $\l$.
  The term $e^{2rs +ts(D\cdot\S)} e^{S_r(\bl,t(D\cdot\S))s}p_r(\l,t)$ 
  is killed by the differential operator 
$$
  P_r=\bs -2r-t(D\cdot \S)- S_r(\bl, t(D\cdot \S)).
$$ 
  Also we may find examples of $(**)$
  with $\deg_{\l} p_r \geq \k$, for any $\k<x_r$.
\end{lem}

\begin{pf}
  Let us suppose $D\cdot \S=1$. Writing $X=X_1\cup_Y A$, where $A$ is 
  a tubular neighbourhood of $\S$, we are under the conditions $(*)$
  for $v=\p^w(X_1,ze^{tD_1})$ with $D=D_1+\D$. Then by corollary~\ref{cor:Sr}
\begin{equation}
  \left\{ \begin{array}{ll}
   (\a-(4r + 2t) - \hS_r(\b+8))v_r=0, (\b+8)^d v_r=0, & \text{$r$ odd}\\
   (\a-(4r\ima-2t)-\hS_r(\b-8))v_r=0, (\b-8)^d v_r=0, \qquad &\text{$r$ even}
  \end{array} \right.
\label{eqn:A}
\end{equation}
  Put $\k$ for the integer such that $(\b+(-1)^{r+1}8)^{\k} v_r \neq 0$ 
  and $(\b+(-1)^{r+1}8)^{\k+1} v_r=0$. By remark~\ref{rem:nilp}
  we can find examples where $\k$ is any integer with $\k <x_r$.
  Using~\eqref{eqn:=}, we have
$$
  \Dws_X(ze^{tD +s\S+\l x})= \la v, \p^w(A,e^{t\D+s\S+\l x}) \ra=\sum_r
  \la v_r, \p^w(A,e^{t\D+s\S+\l x}) \ra=\sum_r D_r,
$$
  where $D_r=\la v_r, \p^w(A,e^{t\D+s\S+\l x}) \ra$. We translate~\eqref{eqn:A}
  into differential equations satisfied by $D_r$. Let
$$
  {\hat P}_r = \left\{ \begin{array}{ll}
  \bs-(2r+t)-\frac12 \hS_r(-4(\bl-2),t), & \text{$r$ odd}\\
  \bs-(2r\ima-t)-\frac12 \hS_r(-4(\bl+2),t), \qquad &\text{$r$ even}
  \end{array} \right.
$$
  Then ${\hat P}_r D_r=0$. For $r$ odd, 
  ${\hat P}_r D_r=\la v_r, \frac12 (\a-(4r+2t) - \hS_r(\b+8,t))
  \p^w(A,e^{t\D+s\S+\l x}) \ra=0$ and analogously for $r$ even. 
  Also $(\bl+(-1)^r2)^{\k} D_r \neq0$
  and $(\bl+(-1)^r2)^{\k+1} D_r = 0$. Therefore
\begin{equation}
 D_r=\left\{ \begin{array}{ll}
   e^{\left(2r+t+{1\over 2}\hS_r(-4(\bl-2),t)\right)s} (e^{2\l}p_r(\l,t))=
     e^{2rs+ts+2\l}e^{{1\over 2}\hS_r(-4\bl,t)s}p_r(\l,t), & \text{$r$ odd} \\
   e^{\left(2r\ima-t+{1\over 2}\hS_r(-4(\bl+2),t)\right)s} 
     (e^{-2\l}p_r(\l,t))= 
     e^{2r\ima s-ts-2\l}e^{{1\over 2}\hS_r(-4\bl,t)s} p_r(\l,t), \quad
   &\text{$r$ even} \end{array} \right. 
\label{eqn:B}
\end{equation} 
with $p_r(u,t) \in \Ct[u]$ a polynomial of degree $\k$. Hence
$$
  D_X^w(z e^{tD+s\S+\l x})= 
  \sum_{r\> \odd} e^{2rs +ts +2\l}e^{{1\over 2}\hS_r(-4\bl,t)s}p_r(\l,t) +
  \sum_{r\> \even}e^{2r \ima s-ts -2\l}e^{{1\over 2}\hS_r(-4\bl,t) s}p_r(\l,t)+
$$
$$
  +\ima^{-d_0}\sum_{r\> \odd}  
   e^{2r\ima s-ts-2\l}e^{{\ima \over 2}\hS_r(4\bl,\ima t)s}p_r(-\l,\ima t)
   +\ima^{-d_0}
  \sum_{r\> \even}
  e^{2rs+ts +2\l} e^{{\ima\over 2}\hS_{-r}(4\bl,\ima t) s} p_{-r}(-\l,\ima t)
$$
  Using~\eqref{eqn:Sr} we get the statement of the lemma. The case of 
  $D\cdot \S$ arbitrary follows readily.
\end{pf}

\begin{thm}
\label{thm:Ar&Br}
  For any $r\in\ZZ$, we have $S_r(u,t)=A_r(u)+B_r(u)t$, with $A_r,B_r\in\Cu$.
\end{thm}

\begin{pf}
  It suffices to prove the statement of the theorem 
  modulo $u^{\k}$, for any $\k<x_r$. By 
  lemma~\ref{lem:inv} there are examples of $4$-manifolds $X$ in the
  situation of $(**)$ such that 
\begin{equation}
  D_X^w(ze^{tD+s\S+\l x})= 
  e^{2\l} \sum_{r=-(g-1)}^{g-1}  
  e^{2rs +ts(D\cdot\S)}e^{S_r(\bl,t(D\cdot\S))s}p_r(\l,t) +\cdots
\label{eqn:abov}
\end{equation}
  with $p_r(u,t)\in \Ct[u]$ of degree at least $\k$ (we do not
  write explicitly the part corresponding to $e^{-2\l}$). Let us argue that
  we may suppose that $z$ does not appear in~\eqref{eqn:abov}. 
  First note that the term
  $e^{S_r(\bl,t(D\cdot\S))s}p_r(\l,t)$ is a polynomial in
  $\l$ and $s$. The general expression~\eqref{eqn:basic}
  for the Donaldson invariants of $X$ gives
\begin{equation}
  D^w_X(e^{t D+s\S+\l x})=e^{Q(tD+s\S)/2+2\l}\sum_K p_K(tD+s\S,\l) 
  e^{K\cdot (tD+s\S)}+ \cdots,
\label{eqn:above}
\end{equation}
  where $p_K(tD+s\S,\l)$ is a polynomial in all $t$, $s$ and $\l$. 
  Comparing~\eqref{eqn:abov} and~\eqref{eqn:above},
  there must be some basic class $K_0$ with $K_0\cdot \S=2r$ such that 
  $\deg p_{K_0} \geq \k$, i.e.\
  $d(K_0)/2 \geq \k$ (recall that we are considering the degree
  with respect to $\l$). Now we choose $D\in H_2(X;\ZZ)$ generic 
  such that $D\cdot \S$ is odd and all values $K \cdot D$ ($K$ 
  basic class of $X$) are different.
  Applying lemma~\ref{lem:inv} to this new $D$ and $z=1$ we have
\begin{equation}
  D_X^w(e^{tD+s\S+\l x})= 
  e^{2\l} \sum_{r=-(g-1)}^{g-1}  
  e^{2rs +ts(D\cdot\S)}e^{S_r(\bl,t(D\cdot\S))s}p_r(\l,t) +\cdots,
\label{eqn:D}
\end{equation}
  for some $p_r(\l,t)$ (different from the previous ones).
  Comparing to~\eqref{eqn:above},
  $e^{S_r(\bl,t(D\cdot\S))s}p_r(\l,t)= \sum_{\{K|K\cdot \S=2r\}}
  e^{Q(tD)/2+K\cdot tD} p_K(tD+s\S,\l)$. By the choice of $D$ there are
  no cancelations in this sum, so the result is a polynomial in $\l$ of 
  degree at least $\k$. Hence $\deg_{\l} p_r\geq \k$.

  Now change $D$ to $D'=aD+b\S$, $a,b \in \ZZ$, $a\neq 0$,
  such that $(D')^2=0$ and 
  either $a$ or $b$ is odd. Then~\eqref{eqn:D} is again satisfied
  for this $D'\in H_2(X;\ZZ)$.
  Moreover if $a$ is odd, then $D'\cdot \S \equiv 1\pmod 2$ 
  so, changing $w$ by $w+\S$ if necessary, we have $w\cdot D'\equiv 1\pmod 2$.
  If $a$ is even, $b$ is odd and then $w\cdot D'\equiv 1\pmod 2$.
  So we may suppose that~\eqref{eqn:D} holds with $D^2=0$ and 
  $w\cdot D\equiv 1\pmod 2$.

  Now~\eqref{eqn:above} becomes
$$
  D^w_X(e^{tD+s\S+\l x})=e^{2\l} \sum_{r,p} 
  e^{(D\cdot \S)ts+2rs+2pt} f_{rp}(t,s,\l)+\cdots,
$$
  for some polynomials $f_{rp} \in \CC[t,s,\l]$. For some integer $p$
  the degree of $f_{rp}$ in the variable $\l$ is at least $\k$. Then 
  lemma~\ref{lem:inv} applied to both embedded surfaces $\S$ and $D$
  implies that $e^{(D\cdot \S)ts+2rs+2pt} f_{rp}(t,s,\l)$ is killed by the 
  differential operators $P_r=\bs -2r-t(D\cdot \S)- 
  S_r(\bl, t(D\cdot \S))$ and
  $Q_p={\bd \over\bd t} -2p-s (D\cdot \S)- S_p(\bl, s(D\cdot \S))$. 
  Henceforth
$$
  [P_r,Q_p] =(D\cdot\S) \left(S_r'(\bl,t(D\cdot \S))-
  S_p'(\bl,s(D\cdot \S))\right)
$$
  kills $f_{rp}(t,s,\l)$ (the prime denotes derivative with respect to
  the second variable). Writing $S_r(u,t)=\sum S_{r,i}(u)t^i$ we see that
  $S_{r,i}(\bl)$ must kill $f_{rp}$ for $i\geq 2$. Therefore
  $S_{r,i}(u)=O(u^{\k+1})$, for $i\geq 2$, as required.
\end{pf}

The proof of theorem~\ref{thm:Ar&Br} also yields that
$x_r= 1+ \sup\{d(K)/2|K$ is a basic class of some 
$4$-manifold $X$ with $b_1=0$ and $b^+>1$ such that there 
is an embedded $\S\subset X$ odd in homology, with 
$K\cdot \S=2r\}$.

\begin{prop}
\label{prop:xr}
  Let $x$ be the supremum of the orders of finite type of
  $4$-manifolds with $b_1=0$ and $b^+>1$. Then 
  $x_r=x$ for all $r\in \ZZ$.
\end{prop}

\begin{pf}
  Fix $r\in \ZZ$. Clearly $x_r \leq x$. Let us see the opposite inequality.
  Let $\k$ be any integer with $\k <x$. Then there exists a 
  $4$-manifold $X$ with $b_1=0$ which has at least one basic class $K$ with
  $d(K)/2\geq \k$. Choose $D\in H_2(X;\ZZ)$ with $D^2>0$
  and take $D'= MD$ where $M$ is a large integer.
  Then blow-up $X$ at $N=(D')^2$ points to get $\tilde X=X \# 
  N \overline{\CP}^2$ with exceptional divisors $E_1,\ldots, E_N$. 
  Consider $\S=D'-E_1-\cdots -E_N$ which has $\S^2=0$
  and it is odd in homology.
  $\tilde K=K +a_1E_1 +\cdots+a_NE_N$, $a_i=\pm 1$, is a basic class of 
  $\tilde X$ and from~\cite{adjuncti} it is $d(\tilde K)/2\geq d(K)/2$. 
  Now $\tilde K \cdot \S= K\cdot D'+a_1+\cdots+ a_N$. 
  Choosing $a_i$ conveniently, we can get any even number between
  $MK\cdot D- M^2D^2$ and $MK\cdot D+ M^2D^2$. With $M$ large enough
  we cover all even numbers.
  Therefore there is a basic class $\tilde K$ for some blow-up such that
  $d(\tilde K)/2\geq d(K)/2\geq \k$ and $\tilde K \cdot \S=2r$.
  So $x_r> d(\tilde K)/2 \geq \k$ and hence $x_r\geq x$. This concludes
  the proof. 
\end{pf}

\begin{rem}
\label{rem:help}
  The construction in the proof of proposition~\ref{prop:xr} can be
  adapted to find a $4$-manifold $X$ with two embedded surfaces
  $\S_i\subset X$, $i=1,2$, odd in homology and with $\S_1^2=0$, 
  $\S_2^2=0$, $\S_1\cdot \S_2=0$, and a basic class $K$ such that 
  $K\cdot \S_1=2r_1$ and $K\cdot \S_2=2r_2$, for any pair of integers 
  $r_1$ and $r_2$. We only need to start with $D_1, D_2\in H_2(X;\ZZ)$ 
  with $D_1^2>0$, $D_2^2>0$ and $D_1\cdot D_2=0$, take large multiples 
  $D_1'=MD_1$, $D_2'=MD_2$ and blow-up at $N=N_1+N_2$ points, 
  with $N_i=(D_i')^2$. Then $\S_1=D'_1-E_1-\cdots -E_{N_1}$ and 
  $\S_2=D'_2-E_{N_1+1}-\cdots E_{N}$ will work.
\end{rem}

\begin{thm}
\label{thm:A&B}
  $A_r(u)=rA_1(u)$ and $B_r(u)=B_0(u)$,
  for any $r\in \ZZ$.
\end{thm}

\begin{pf}
  It suffices to prove the statement of the theorem 
  modulo $u^{\k}$ for any $\k<x$. There are two cases
\begin{itemize}
  \item $r$ odd. By lemma~\ref{lem:inv} and theorem~\ref{thm:Ar&Br} 
  we may find $X$ satisfying $(**)$ with $z=1$ and
$$
  D_X^w(e^{tD+s\S+\l x})= 
  e^{2\l} \sum_{r=-(g-1)}^{g-1}  
  e^{2rs +ts(D\cdot\S)} e^{A_r(\bl)s+B_r(\bl)ts(D\cdot\S)}
  p_r(\l,t) +\cdots,
$$
  $\deg_{\l} p_r \geq \k$. Let $\S'=n\S$, for $n$ odd,
  and $D'=D/n$. Then $\S'$ is represented by an embedded surface of genus
  $gn$ with $(\S')^2=0$ and $w\cdot \S'\equiv 1\pmod 2$. Then
$$
  D_X^w(e^{tD'+s\S'+\l x})= 
  e^{2\l} \sum_{r=-(g-1)}^{g-1}  
  e^{2rns +ts(D\cdot\S)}e^{A_r(\bl)ns+B_r(\bl)ts(D\cdot\S)}
  p_r(\l,t/n) +\cdots
$$
  By lemma~\ref{lem:inv} applied to $\S'$,
  $e^{2rns +ts(D\cdot\S)}e^{A_r(\bl)ns+B_r(\bl)ts(D\cdot\S)} p_r(\l,t/n)$
  is killed by $P_{rn}$. So it must be
  $A_{rn}=nA_r$ and $B_{rn}=B_r$ (for $n$ odd) modulo $u^{\k}$.
  This proves that $A_r=rA_1$ and $B_r=B_1$ for $r$ odd. Also $A_0=0$.

\item $r$ even. By remark~\ref{rem:help}, there is 
  $4$-manifold $X$ with a basic 
  class $K$ with $d(K)/2 \geq \k$ such that there exist
  $\S,D\in H_2(X;\ZZ)$, $w\in H^2(X;\ZZ)$ with $\S^2=0$, 
  $D^2=0$, $\S \cdot D=0$, $w\cdot \S \equiv 1\pmod 2$ and
  $w\cdot D \equiv 1\pmod 2$, satisfying
  $K\cdot \S=2r$ and $K\cdot D=2$. For $H \in H_2(X)$ 
$$
  D^w_X(e^{s\S+l D+t H+\l x})= e^{2\l}\sum_{r,p}
  e^{2rs+2pl+(s\S+l D)\cdot tH} f_{rp}(s,l,t,\l)+\cdots,
$$
  where $f_{rp}(s,l,t,\l)$ is polynomic in $s$, $l$ and $\l$. 
  Choosing $H$ generic, the term $f_{r1}(s,l,t,\l)$ has 
  degree bigger or equal than $\k$ in the variable $\l$.
  By lemma~\ref{lem:inv} applied to $\S$, $D$ and $\S'=\S-rD$ 
  (using here that $r$ is even), 
  $e^{2rs+2l+(s\S+l D)\cdot tH} f_{r1}(s,l,t,\l)$
  is a solution to 
\begin{eqnarray*}
  & & \bs -2r-t(H \cdot \S)- A_r(\bl)-B_r(\bl)t(H\cdot \S), \\
  & & {\bd\over \bd l}-2-t(H\cdot D)-A_1(\bl)-B_1(\bl)t(H\cdot D)   \quad \text{and}\\
  & & \bs-r{\bd\over \bd l}  -t(H\cdot \S')-B_0(\bl)t(H\cdot \S').
\end{eqnarray*}
  Hence
$$
  \left(A_r(\bl)-rA_1(\bl)\right) -
  \left(B_r(\bl)-B_0(\bl)\right)t(H\cdot \S) -
  \left(B_1(\bl)-B_0(\bl)\right) rt(H\cdot D)
$$
  kills $f_{r1}(s,l,t,\l)$. As this works for $H$ generic, all three 
  expressions in between parentheses
  vanish, i.e.\ $A_r=rA_1$ and $B_r=B_1=B_0$.
\end{itemize}
\end{pf}
  
{\em Proof of theorem~\ref{thm:main2}.\/}

We put $A(u)=1+\frac12 A_1(u)$, $B(u)=1+B_0(u)$. 
Then the expression in lemma~\ref{lem:inv} becomes
$$
  D_X^w(ze^{tD+s\S+\l x})= e^{2\l}\sum_{r=-(g-1)}^{g-1}  
   e^{A(\bl)2rs +B(\bl)ts(D\cdot\S)} p_r(\l,t)+\cdots,
$$
and $P_r=\bs-A(\bl)2r-B(\bl)t(D\cdot\S)$. 

As in the proof of theorem 6 in~\cite{basic}, it is enough to prove the
statement for some blow-up $\tilde X= X\# m\overline{\CP}^2$ 
of $X$ with exceptional divisors
$\seq{E}{1}{m}$ and some $\tilde w \in H^2(\tilde X;\ZZ)$ of the form
$\tilde w=w +\sum a_i E_i$, $a_i\in \ZZ$. 
So by step 2 of proof of theorem 6 in~\cite{basic}, it is enough to
prove theorem~\ref{thm:main2} in the case where 
$X$ is a $4$-manifold with $b_1=0$, $b^+>1$ such that there is 
$w \in H^2(X;\ZZ)$ and a subgroup $H=\la\seq{\S}{1}{n}\ra \subset
\bar H_2(X;\ZZ)$ with $2\bar H_2(X;\ZZ)\subset H$, $\S_j^2=0$
and $w \cdot \S_j \equiv 1 \pmod 2$, $1 \leq j\leq n$. Represent
$\S_j$ by an embedded surface of genus $g_j$ and apply
lemma~\ref{lem:inv} simultaneously to all $\S_j$ to get
$$
  D^w_X(e^{t_1\S_1 +\cdots+ t_n \S_n+\l x})=e^{2\l} 
 \hspace{-5mm} \sum_{-(g_j-1)\leq r_j \leq g_j-1\atop 1\leq j\leq n}
  \hspace{-5mm}
  e^{A(\bl)(2r_1t_1 +\cdots+2r_nt_n) +B(\bl)Q(t_1\S_1 +\cdots+ t_n
  \S_n)/2} P_{r_1\ldots r_n,w}(\l)+\cdots
$$
where $P_{r_1\ldots r_n,w}$ is a polynomial in $\l$. This yields the
statement.
\hfill $\Box$

\begin{rem}
\label{rem:x}
  Let $x$ denote the supremum of the orders of finite type of 
  $4$-manifolds with $b_1=0$ and $b^+>1$ (so that the simple type
  conjecture reads as $x=1$). If $x<\infty$ then the series 
  $A(u)$ and $B(u)$ appearing in theorem~\ref{thm:main2} 
  are uniquely determined modulo $u^x$. The terms of degree higher 
  than $x-1$ are irrelevant for the Donaldson invariants of $4$-manifolds 
  with $b_1=0$, $b^+>1$, since $\deg P_K \leq x-1$ for
  any basic class $K$.
  If $x=\infty$ then $A(u)$ and $B(u)$ are uniquely determined.
\end{rem}

\section{Review of modular forms}
\label{sec:4}

We are going to review the notations and results on modular
forms that we shall use in this paper. For general background 
on modular forms and elliptic functions 
see the books~\cite{Ch}\cite{Lawden}.
In general, we follow the notations of~\cite{GZ} but we 
shall point out the differences.

Let $\HH=\{\tau\in \CC | \text{Im} (\tau)>0\}$ be the complex upper
half-plane. For $\tau\in \HH$ let $q=e^{2\pi i \tau}$ and 
$q^{1/n}=e^{2\pi i \tau/n}$. 
We denote $\s_1(n)=\sum_{d|n}d$ and by 
$\s_{1}^{\odd}(n)$ the sum of the odd divisors of $n$.
Let
$$
  G_2(\tau)=-\frac{1}{24}+\sum_{n>0}\sigma_{1}(n)q^n
$$
be the classical Eisenstein series. This is a quasi-modular form,
i.e.\
$$
G_2(\frac{a\tau+b}{c\tau+d})=(c\tau+d)^2G_2(\tau)- \frac{c
(c\tau+d)}{4\pi\ima},
$$
for $\left(\begin{array}{cc}a&b\\ c&d\end{array}\right)
 \in SL(2,\ZZ)$. Let
\begin{align*}
 \theta(\tau,z) &=-\ima \sum_{n=-\infty}^{\infty}(-1)^nq^{{1\over 2}(n+{1\over 
 2})^2} e^{2\pi\ima (n+{1\over 2})z} \\
\theta_1(\tau,z) &=\sum_{n=-\infty}^{\infty} q^{{1\over 2}(n+{1\over 
 2})^2} e^{2\pi\ima (n+{1\over 2})z} \\
\theta_2(\tau,z) &=\sum_{n=-\infty}^{\infty} (-1)^n q^{{1\over 2}n^2} 
  e^{2\pi\ima nz} \\
\theta_3(\tau,z) &=\sum_{n=-\infty}^{\infty} q^{{1\over 2}n^2} 
  e^{2\pi\ima nz}
\end{align*}
be the theta functions. In general we shall also write
$\theta_i(\tau)=\theta_i (\tau,0)$, $i=1,2,3$.
These are the theta functions of~\cite{Ch}
and coincide with those of~\cite{GZ} except for the factor of
$-\ima$ in $\h$. The theta functions are labelled as $\h_1$,
$\h_2$, $\h_4$, $\h_3$ in~\cite{Lawden} and~\cite{Moore-Witten}
(moreover in~\cite{Lawden} the factor of $\pi$ in the
exponentials does not appear, so that e.g.\
$\h'(\tau,0)$ for us corresponds to $\pi\h'(\tau,0)$
in~\cite{Lawden}).
Also let $\eta(\tau)=q^{1/24}\prod_{n>0}(1-q^n)$ be the Dirichlet
eta function, so $\h'(\tau,0)=2\pi \eta(\tau)^3=\pi
\h_1(\tau)\h_2(\tau)\h_3(\tau)$. By~\cite[Chapter V, 
theorem 5]{Ch}, it is
$\h'''(\tau,0)=4\pi\ima {\bd\over\bd\tau}\h'(\tau,0)$ whence
\begin{equation}
 \frac{\h'''(\tau,0)}{\h'(\tau,0)}= 4\pi\ima {d\over d\tau}
  \log \h'(\tau,0)= 12\pi\ima {d\over d\tau}
  \log \eta(\tau)= -24\pi^2 q {d\over d q} \log \eta(\tau)=
  24 \pi^2 G_2(\tau),
\label{eqn:lp}
\end{equation}
where $G_2(\tau)=-q {d\over d q} \log \eta(\tau)$
follows by term by term logarithmic
differentiation in the definition of $\eta(\tau)$.

Finally let $e_1$, $e_2$, $e_3$ be the $2$-division values 
of the Weierstrass 
function at $1/2$, $\tau/2$ and $(1+\tau)/2$ respectively, i.e.\
\begin{align*}
 e_1(\tau) &=-\frac{1}{6}-4\sum_{n>0}\sigma^{\text{odd}}_1(n)q^{n} 
  =-\frac{1}{12}(\theta_2^4(\tau) +\theta_3^4(\tau)) \\
 e_2(\tau) &=\frac{1}{12}+2\sum_{n>0}\sigma^{\text{odd}}_1(n)q^{n/2}
  =\frac{1}{12}(\theta_1^4(\tau) +\theta_3^4(\tau)) \\
 e_3(\tau) &=\frac{1}{12}+2\sum_{n>0}(-1)^n\sigma^{\text{odd}}_1(n)q^{n/2} 
  =-\frac{1}{12}(\theta_1^4(\tau) - \theta_2^4(\tau)) 
\end{align*}
The right hand side equalities follow from~\cite[Chapter V,
 (5.6)]{Ch}, taking into account that our functions
$e_i(\tau)$ are those of~\cite{Ch} divided by a factor $-4\pi^2$.
Therefore by~\cite[Chapter V, (5.2)]{Ch}
\begin{equation}
  4\pi\ima {d\over d\tau} \log \h_i(\tau)=
 \frac{\h''_i(\tau,0)}{\h_i(\tau,0)}= 
 8 \pi^2 G_2(\tau) +  4\pi^2e_i(\tau), \qquad i=1,2,3.
\label{eqn:lpi}
\end{equation}
  
In~\cite{GZ} the following functions are introduced
\begin{align}
  f(\tau) &=e^{\pi\ima 3/4}\frac{\eta(\tau)^3}{\h_3(\tau)}=
  {1\over 2}e^{\pi\ima 3/4}\h_1(\tau)\h_2(\tau), \nonumber\\
  U(\tau)&=-\frac{3e_3(\tau)}{f(\tau)^2}, \nonumber\\
  G(\tau)&= \frac{G_2(\tau)+e_3(\tau)/2}{f(\tau)^2}.
\label{eqn:magnetic}
\end{align}

\begin{thm}[\cite{FS-bl}]
\label{thm:blow-up}
  Let $B(x,t)=e^{-\frac{t^2x}{6}} \sigma_3(t)$ and $S(x,t)=
  e^{-\frac{t^2x}{6}} \sigma(t)$, where $\sigma$ and $\sigma_3$ are
  $\sigma$-functions associated to the elliptic curve with 
  $g_2=4(\frac{x^2}{3}-1)$ and $g_3=\frac{8x^3-36x}{27}$.
  For any $4$-manifold $X$ with $b^+>1$ and $b_1=0$, let $\tilde X=
  X\# \overline{\CP}^2$ denote its blow-up and let $E$ stand for
  the exceptional divisor. For $w\in H^2(X;\ZZ)$ we have
\begin{align*}
  & D^w_{\tilde X}(e^{tE+sD+\l x})= D^w_X(B(x,t)e^{sD+\l x}), \\
  & D^{w+E}_{\tilde X}(e^{tE+sD+\l x})= D^w_X(S(x,t)e^{sD+\l x}),
\end{align*}
  for all $D \in H_2(X)$. \hfill $\Box$
\end{thm}

\begin{prop} {\em 
(\cite[proposition 2.19]{GZ}\cite[section 6]{Moore-Witten})}
\label{prop:blow-up}
We have
\begin{align*}
  B(U(\tau),t)=& e^{t^2 G(\tau)}
     \frac{\theta_3(\tau,\frac{t}{2\pi\ima f(\tau)})}{\theta_3(\tau,0)}, \\
  S(U(\tau),t)=& e^{-3\pi \ima/4} e^{t^2 G(\tau)}
     \frac{\theta(\tau,\frac{t}{2\pi\ima f(\tau)})}{\theta_3(\tau,0)}.
\end{align*}
\end{prop}

\begin{pf}
  This is an exercise in modular forms. We shall prove only the second 
  equality, the first one being similar.  
  Put $\omega_1=2\pi\ima f(\tau)$ and $\omega_2=2\pi\ima \tau f(\tau)$
  and let us find the invariants $g_2$ and $g_3$ associated to the lattice
  $L=\ZZ\o_1+\ZZ\o_2$. For brevity, set 
  $x_i=\h_i^4$, $i=1,2,3$, so that $x_3=x_1+x_2$. Then $12e_3=x_2-x_1$
  and $16 f^4=-x_1x_2$. Using~\cite[Chapter 6, exercise 4]{Lawden}
  (and noting that in~\cite{Lawden} the basic periods are
  $2\o_1$, $2\o_3$), we have 
\begin{align*}
  g_2(L) &= {2\pi^4\over 3\o_1^4}(\h_1^8+\h_2^8+\h_3^8)= 
   {1\over 24 f^4}(\h_1^8+\h_2^8+\h_3^8) = 
   -2 {x_1^2+x_2^2+x_3^2 \over 3 x_1x_2} = \\
  &= -4 {x_1^2+x_2^2+x_1x_2 \over 3 x_1x_2} = 
   -4 \left( {(12e_3)^2 \over 3 x_1x_2} + 1 \right)= 
   4 \left( {(-3e_3/ f^2)^2\over 3} -1\right) = 4({U^2\over 3}-1), \\
  g_3(L) & = {4 \pi^6\over 27\o_1^6}
  (\h_1^4+\h_3^4)(\h_2^4+\h_3^4)(\h_2^4-\h_1^4)=
   {4 (x_1+x_3)(x_2+x_3) \over 27x_1x_2} {12 e_3 \over 4f^2} =\\
  &=- {4 (2x_1+x_2)(2x_2+x_1) \over 27x_1x_2} U = {1\over 27} 
  \left( { 8(x_1-x_2)^2\over -x_1x_2} -36 \right) U= {1\over 27}(8U^2-36) U.
\end{align*}
Therefore putting $x=U(\tau)$ the $\s$-functions of 
theorem~\ref{thm:blow-up} correspond to the lattice $L$. Now
by~\cite[Chapter V, theorem 2]{Ch},
$$
  \sigma(t)= {\omega_1\over \theta'(\tau,0)}e^{\eta_1 t^2/\o_1}
  \theta(\tau, {t \over\omega_1})  =
  {2\pi\ima f(\tau)\over \pi \h_1(\tau)\h_2(\tau)\h_3(\tau)}
  e^{\eta_1 t^2/\o_1} \theta(\tau,{t\over 2\pi\ima f(\tau)}) =
  \ima e^{3\pi \ima/4} e^{\eta_1 t^2/\o_1} 
  {\h (\tau,\frac{t}{2\pi\ima f(\tau)}) \over \h_3(\tau)}
$$
and by~\cite[formula (6.2.7)]{Lawden} and~\eqref{eqn:lp},
$\eta_1=-{1\over 6\o_1}{\h'''(\tau,0)\over \h'(\tau,0)}=
 -{4 \pi^2 \over \o_1} G_2(\tau)$. So
$$
  S(U(\tau),t) = e^{-t^2U(\tau)/6} \sigma(t)=
  exp({e_3(\tau)\over 2f(\tau)^2}t^2 + 
  {G_2(\tau)\over f(\tau)^2} t^2) 
   e^{-3\pi \ima/4} {\h (\tau,\frac{t}{2\pi\ima f(\tau)}) 
   \over \h_3(\tau,0)}= e^{-3\pi \ima/4} e^{t^2 G(\tau)}
  \frac{\h(\tau,\frac{t}{2\pi\ima f(\tau)})}{\theta_3(\tau,0)}.
$$
\end{pf}

\section{Proof of Theorem \protect\ref{thm:main}}
\label{sec:5}

In this section we shall prove our main theorem. For this
let us rewrite the expression of theorem~\ref{thm:main2} 
as a coefficient of $q^0$ of a Laurent power series in 
$q=e^{2\pi\ima\tau}$.

\begin{lem}
\label{lem:mod}
  Let $V(\tau)=\sum_{n\geq 0} a_nq^n$ be a power series in $q$
  with $a_0=2$ and $a_1 \neq 0$. 
  Let $N(\tau)=A(V(\tau)-2)$ and $M(\tau)=B(V(\tau)-2)$ which are 
  formal power series in $q$. Then for any $4$-manifold $X$ with 
  $b^+>1$ and $b_1=0$, and $K\in H^2(X;\ZZ)$ a basic class for $X$,
  there exists a uniquely determined polynomial $f_K(\tau)$ in $q$ 
  of degree $d(K)/2$, such that
  $$
  e^{2\l} e^{B(\bl)Q(D)/2+
  A(\bl) K\cdot D} P_K(\l) = \left[q^{-d(K)/2} e^{V(\tau)\l +
  M(\tau)Q(D)/2+N(\tau) K\cdot D} f_K(\tau) \right]_{q^0},
  $$
  for any $D\in H_2(X)$, where $[F]_{q^0}$ stands for the 
  coefficient of $q^0$ of a Laurent power series $F(q)$.
\end{lem}

\begin{pf}
  Let $P(D,\l)$ stand for the left hand side. Then $P$ is 
  characterised by being a solution of
\begin{equation} 
  \left\{
\begin{array}{l}
  \left(\bt -(D\cdot H) B(\bl-2) -(K\cdot H)A(\bl-2)\right)
  \bigg|_{t=0} P(D+t H,\l)=0, \\
  \left(\bl-2\right)^{\frac{d(K)}{2}+1} P(D,\l)=0,
\end{array} \right.
\label{eqn:sol}
\end{equation}
  for any $D, H \in H_2(X)$. Now
$$ 
  \left(\bt -(D\cdot H) B(\bl-2)-(K\cdot H)A(\bl-2)\right)\bigg|_{t=0}
  \left[q^{-d(K)/2} e^{V(\tau)\l +
  M(\tau)Q(D+tH)/2+N(\tau) K\cdot (D+tH) } f_K(\tau) \right]_{q^0}=
$$
$$
  = \bigg[q^{-d(K)/2} \left( (D\cdot H)M(\tau)+(K\cdot H) N(\tau)
  -(D\cdot H) B(V(\tau)-2) -(K\cdot H)A(V(\tau)-2)\right) \cdot
   \qquad
$$
$$
  \qquad \cdot\, e^{V(\tau)\l +
  M(\tau)Q(D)/2+N(\tau)K\cdot D} f_K(\tau) \bigg]_{q^0}=0,
$$
  and
$$
  \left(\bl-2\right)^{\frac{d(K)}{2}+1} 
  \left[q^{-d(K)/2} e^{V(\tau)\l +
  M(\tau)Q(D)/2+N(\tau) K\cdot D } f_K(\tau) \right]_{q^0}=
$$
$$
  =\left[\left(V(\tau)-2\right)^{\frac{d(K)}{2}+1} 
  q^{-d(K)/2} e^{V(\tau)\l +
  M(\tau)Q(D)/2+N(\tau) K\cdot D} f_K(\tau) \right]_{q^0}=0,
$$
  since the series $V(\tau)-2$ equals 
  $a_1q$ plus higher order terms. 
  The condition $a_1\neq 0$ ensures that varying $f_K(\tau)$ in 
  the polynomials of degree $d(K)/2$ we get a basis for the 
  solutions to~\eqref{eqn:sol}.
  So there exists exactly one $f_K(\tau)$ as required. 
\end{pf}

\begin{rem}
\label{rem:l}
  Although in lemma~\ref{lem:mod} we are supposing that 
  $f_K(\tau)$ is a polynomial in $q$, we may let $f_K(\tau)$
  be a power series instead. The coefficients will be
  determined only up to the power $d(K)/2$.
\end{rem}

To find $V$, $N$ and $M$ we shall make use of their universality
and the blow-up formula. For this we need to re-express 
the blow-up formula. Let us define the following functions 
similar to~\eqref{eqn:magnetic}
\begin{align}
  h(\tau)& =\frac12 \h_2(\tau)\h_3(\tau) =\frac12 -2q+2q^2+\cdots
  \nonumber\\
  V(\tau)& =\frac{-3e_1(\tau)}{h(\tau)^2}=2+64q+512q^2+\cdots, 
  \nonumber\\
  T(\tau)& =-\frac{G_2(\tau)+e_1(\tau)/2}{h(\tau)^2}=\frac12+8q+30q^2+\cdots
\label{eqn:electric}
\end{align}

\begin{lem}
\label{lem:V}
We have
\begin{align*}
  B(V(\tau),t) &= e^{-t^2 T(\tau)}
  \frac{\h_1(\tau,\frac{t}{2\pi\ima h(\tau)})}{\h_1(\tau,0)} \\
  S(V(\tau),t) &= \ima\,e^{-t^2 T(\tau)}
  \frac{\h(\tau,\frac{t}{2\pi\ima h(\tau)})}{\h_1(\tau,0)}
\end{align*}
\end{lem}

\begin{pf}
We consider the transformation $\tau'= 1-{1\over\tau}$.
This is what in physical terms corresponds to changing from
``magnetic'' variables to ``electric'' variables. Using the
modular behaviour of the theta functions~\cite[Chapter V,
section 8]{Ch} and the
quasi-modular behaviour of $G_2$, we have
\begin{align*}
 & f(\tau')={\ima\over 2} \tau\h_2(\tau)\h_3(\tau)=\ima\tau h(\tau),\\
 & G(\tau')=-\frac{G_2(\tau)+e_1(\tau)/2}{h(\tau)^2} 
  +\frac{1}{4\pi\ima \tau h(\tau)^2}, \\
 & U(\tau')=\frac{3e_1(\tau)}{h(\tau)^2}=-V(\tau).
\end{align*}
Applying the transformation to the formula in
proposition~\ref{prop:blow-up} we get
$$
  B(-V(\tau), t) = 
  exp(-t^2 \frac{G_2(\tau)+e_1(\tau)/2}{h(\tau)^2} + 
  \frac{t^2}{4\pi\ima \tau h(\tau)^2})
  \frac{\sqrt{\frac{\tau}{\ima}}e^{\frac{\pi\ima}{\tau}\left(
  \frac{t}{-2\pi h(\tau)}\right)^2}
  \theta_1(\tau,\frac{t}{-2 \pi h(\tau)})}
  {\sqrt{\frac{\tau}{\ima}} \theta_1(\tau,0)},
$$
and simplifying 
$$
  B(-V(\tau),t)= e^{t^2T(\tau)}
   \frac{\theta_1(\tau,\frac{t}{-2 \pi h(\tau)})}{\h_1(\tau,0)}.
$$
Using that $B(V(\tau),t)=B(-V(\tau),\ima t)$ we get the
statement. The second line is proved analogously.
\end{pf}

\begin{prop}
\label{prop:V}
  Let $V(\tau)=\frac{-3e_1(\tau)}{h(\tau)^2}$ as above. 
  Then $M(\tau)=2T(\tau)$ and $N(\tau)=1/2h(\tau)$. So 
  for any $4$-manifold $X$ with $b_1=0$ and $b^+>1$ and odd,
  and $w \in H^2(X;\ZZ)$,
  there exist polynomials $f_{K,w}(\tau)$ in $q$ of degree $d(K)/2$,
  for every basic class $K\in H^2(X;\ZZ)$, such that
 $$
  D^w_X(e^{D+\l x})= \left[q^{-d(K)/2} e^{V(\tau)\l +
  T(\tau)Q(D)+ K\cdot D/2h(\tau)} f_{K,w}(\tau) \right]_{q^0} +
$$
$$
  + \ima^{-d_0} \left[q^{-d(K)/2} e^{-V(\tau)\l
  -T(\tau)Q(D)/2+ \ima K\cdot D/h(\tau)} f_{K,w}(\tau) \right]_{q^0},
 $$
  for $D\in H_2(X)$.
\end{prop}

\begin{pf}
  Fix $w\in H^2(X;\ZZ)$. 
  Let $\tilde X=X\# \overline{\CP}^2$ be the blow-up of $X$, 
  with exceptional divisor $E$. Then theorem~\ref{thm:blow-up}
  and lemma~\ref{lem:mod} say that
$$
  D^w_{\tilde X}(e^{tE+D+\l x})= D^w_X(B(x,t)e^{D+\l x})=
  \left[q^{-d(K)/2} B(V(\tau),t) e^{V(\tau)\l +
  M(\tau)Q(D)/2+N(\tau) K\cdot D} f_K(\tau) \right]_{q^0}+\cdots
$$
Now by lemma~\ref{lem:V},
\begin{align*}
 B(V(\tau),t) & =e^{-t^2 T(\tau)}\frac{\h_1(\tau,
  \frac{t}{2\pi\ima h(\tau)})}{\h_1(\tau,0)}= e^{-t^2 T(\tau)}{1
  \over \h_1(\tau)} \sum_{n=-\infty}^{\infty}
  q^{1/8} q^{{1\over2}n(n+1)} e^{(2n+1)t/2h(\tau)}= \\
 &= e^{-t^2 T(\tau)}{1 \over  \tilde\h_1(\tau)} 
  \sum_{n=0}^{\infty} q^{{1\over2}n(n+1)} (e^{(2n+1)t/2h(\tau)}
  +e^{-(2n+1)t/2h(\tau)}),
\end{align*}
where $\tilde\h_1(\tau)=q^{-1/8}\h_1(\tau)=2+2q+2q^3+\cdots$, 
which is a power series in $q$. Substituting
into the expression above
$$
  D^w_{\tilde X}(e^{tE+D+\l x})=
  \sum_K \sum_{n=0}^{\infty} 
  \bigg[q^{-d(K)/2+{1\over2}n(n+1)} e^{V(\tau)\l}
  e^{-T(\tau)t^2+ M(\tau)Q(D)/2} \cdot\qquad
$$
$$
  \qquad\cdot\left( e^{(2n+1)t/2h(\tau)+ N(\tau) K\cdot D} +
  e^{-(2n+1)t/2h(\tau)+ N(\tau) K\cdot D} \right)
  {f_K(\tau) \over \tilde\h_1(\tau)} \bigg]_{q^0}+\cdots
$$

As $Q(tE+D)=-t^2+Q(D)$, it must be $M(\tau)=2 T(\tau)$.
The expression inside the summatory is non-vanishing if
and only if $n(n+1)\leq d(K)$. So we actually have a
finite sum. From~\cite{basic} the basic classes of $\tilde X$ 
are of the form $\tilde K=K + (2m+1)E$,
where $K$ is a basic class of $X$ and $m\in\ZZ$. The expression in
the summatory correspond to the pair of basic classes
$\tilde K=K \pm (2n+1)E$, since in this case
$\tilde K\cdot (tE+D)=\mp (2n+1)t + K\cdot D$.
So it is $N(\tau)=1/2h(\tau)$. 
Note also that $\tilde K=K\pm (2n+1)E$
is a basic class for $\tilde X$ if and only if $K$ is a basic
class for $X$ and $d(\tilde K)=d(K)-n(n+1)\geq 0$. 
\end{pf}

To finish the proof of theorem~\ref{thm:main} it only remains
to prove that

\begin{prop}
\label{prop:w}
  $f_{K,w}(\tau)=(-1)^{K\cdot w +w^2\over 2}f_K(\tau)$, for
  any $w \in H^2(X;\ZZ)$. Here $f_K(\tau)=f_{K,0}(\tau)$.
\end{prop}

\begin{pf}
Let $\tilde X=X\# \overline{\CP}^2$ be the blow-up of $X$, 
with exceptional divisor $E$ and let $w\in H^2(X;\ZZ)$.
We are going to relate $f_{K,w}$ with $f_{K\pm E,w}$
and $f_{K\pm E,w+E}$, for a basic class $K$ of $X$. 
From the proof of proposition~\ref{prop:V} we see
that $f_{K+E,w}(\tau)=f_{K-E,w}(\tau)=f_K(\tau)/
\tilde\h_1(\tau)$. On the other hand
$D^{w+E}_{\tilde X}(e^{tE+D+\l x})= D^w_X(S(x,t)e^{D+\l x})$
and
\begin{align*}
 S(V(\tau),t) & = \ima \,e^{-t^2 T(\tau)}\frac{\h(\tau,
  \frac{t}{2\pi\ima h(\tau)})}{\h_1(\tau,0)}=
  e^{-t^2 T(\tau)}{1
  \over \h_1(\tau)} \sum_{n=-\infty}^{\infty} (-1)^n
  q^{1/8} q^{{1\over2}n(n+1)} e^{(2n+1)t/2h(\tau)}= \\
 &= e^{-t^2 T(\tau)}{1 \over \tilde\h_1(\tau)} 
  \sum_{n=0}^{\infty} (-1)^n q^{{1\over2}n(n+1)} 
  (e^{(2n+1)t/2h(\tau)} - e^{-(2n+1)t/2h(\tau)}),
\end{align*}
so $f_{K-E,w+E}(\tau)=f_K(\tau)/\tilde\h_1(\tau)=
-f_{K+E,w+E}(\tau)$. Therefore $f_{K+E,w+E}=-f_{K+E,w}$
and $f_{K-E,w+E}=f_{K-E,w}$.
The arguments in step 4 of the proof of theorem 2 in~\cite{basic}
now carry over.
\end{pf}

\section{Proof of theorem~\protect\ref{thm:main3}}
\label{sec:6}

By using lemma~\ref{lem:mod} in the expression of 
theorem~\ref{thm:main}, we get that for any $4$-manifold
$X$ with $b_1=0$, $b^+>1$ and odd, there are polynomials
$P_K(\l)$ of degree $d(K)/2$ such that 
for any $w\in H^2(X;\ZZ)$, we have 
  $$
  D^w_X(e^{tD+\l x})= e^{2\l}\sum_K e^{B(\bl)Q(tD)/2+
  A(\bl) K\cdot tD} (-1)^{K\cdot w+w^2\over 2}
  P_K(\l) +\qquad
  $$
  $$
  \qquad +\ima^{-d_0} e^{-2\l}\sum_K 
  e^{-B(-\bl)Q(tD)/2+A(-\bl)\ima K\cdot tD}
  (-1)^{K\cdot w+w^2\over 2} P_K(-\l),
  $$
where $A(u), B(u)\in \Cu$ are universal series satisfying
\begin{equation}
\left\{ \begin{array}{l}
   A(V(\tau)-2)=N(\tau)=1/2h(\tau), \\
   B(V(\tau)-2)=M(\tau)=2 T(\tau).
\end{array}\right.
\label{eqn:98}
\end{equation}

In order to solve~\eqref{eqn:98} we need to review some 
elementary facts on elliptic integrals
as presented in~\cite[Chapters 2 and 3]{Lawden}. The (complete) 
elliptic integrals of first and second kind are
\begin{align*}
  K &= \int_0^1 \frac{dt}{\sqrt{(1-t^2)(1-k^2t^2)}} = \frac12
  \pi \left[ 1+\left(\frac12 \right)^2 k^2 + \left( 
  \frac{1\cdot3}{2\cdot4}\right)^2k^4 +\left(
  \frac{1\cdot3\cdot5}{2\cdot4\cdot6}\right)^2k^6+\cdots\right],\\
  E &= \int_0^1 \sqrt{\frac{1-k^2t^2}{1-t^2}} dt = \frac12
  \pi \left[ 1-\left(\frac12 \right)^2 k^2 -\frac13 \left( 
  \frac{1\cdot3}{2\cdot4}\right)^2k^4 -\frac15 \left(
  \frac{1\cdot3\cdot5}{2\cdot4\cdot6}\right)^2k^6-\cdots\right],
\end{align*}
respectively. These are analytic functions of $k^2$. Here $k$
is termed the modulus of the elliptic integral. The relation
with modular forms come from the classical resolution given by
(see section 2 and formula (3.3.5) in~\cite{Lawden})
\begin{align*}
  k &= k(\tau)=\frac{\h_1^2(\tau)}{\h_3^2(\tau)}, \\
  K &= \frac12 \pi \h_3^2(\tau), \\
  E &= \left(1-
  \frac{\h_2''(\tau,0)}{\pi^2\h_2(\tau,0)\h_3^4(\tau)}\right) K.
\end{align*}
Using~\eqref{eqn:lpi} we get
$$
  4(K-E)K=\frac{\h_2''(\tau,0)}{\h_2(\tau,0)}=8\pi^2 G_2(\tau)+
  4\pi^2 e_2(\tau).
$$
As $e_2-e_1=\frac14 \h^4_3(\tau)=K^2/\pi^2$, we also have
$$
  -4EK=8\pi^2 G_2(\tau)+ 4\pi^2 e_1(\tau).
$$

Now recall the duplication formulae~\cite[section 1.8]{Lawden}
\begin{align*}
  &\h_2^2(2\tau) = \h_2(\tau)\h_3(\tau), \\
  & e_2(2\tau) = - e_1(\tau)/2, \\
  & \h_1(2\tau)\h_3(2\tau) = \h_1^2(\tau)/2
\end{align*}
(the second line follows from the definition).

Put $\tau'=2\tau+1$ and let $k=k(\tau')$ be the modulus 
corresponding to $\tau'$. Then we have
\begin{align}
 h(\tau) &=\frac12 \h_2(\tau)\h_3(\tau)=
  \frac12 \h_2^2(2\tau)=\frac12 \h_3^2(\tau')= K/\pi, \nonumber\\
 V(\tau) & =\frac{-3e_1(\tau)}{h(\tau)^2}=
  \frac{24 e_2(2\tau)}{\h_2^2(2\tau)}= 
  \frac{24 e_3(\tau')}{\h_3^2(\tau')}= 
  -2 \frac{\h_1^4(\tau')-\h_2^4(\tau')}{\h_3^2(\tau')}=
  -4\frac{\h_1^4(\tau')}{\h_3^2(\tau')}+2=-4k^2+2, \nonumber\\
 T(\tau) & =-\frac{G_2(\tau)+e_1(\tau)/2}{h(\tau)^2}=
  -\pi^2 \frac{G_2(\tau')+e_1(\tau')/2 + G_2(\tau')+
  e_2(\tau')/2}{K^2}=-\frac{4(K-E)K-4EK}{8K^2}=\nonumber\\
  & = \frac{2E-K}{2K}.
\label{eqn:2E-K}
\end{align}
In the third line we have used~\eqref{eqn:lpi} to get
$$
G_2(\tau)+\frac{e_1(\tau)}2= 4\pi\ima {d\over d\tau} \log 
\h_1(\tau)= 2\pi\ima \left({d\over d\tau} \log \h_1(2\tau)+
{d\over d\tau} \log \h_3(2\tau) \right)=
$$
$$
=4\pi\ima \left({d\over d\tau'} \log \h_1(\tau')+
{d\over d\tau'} \log \h_2(\tau') \right)=
G_2(\tau')+\frac{e_1(\tau')}2 +
G_2(\tau')+\frac{e_2(\tau')}2.
$$

Using~\eqref{eqn:2E-K} we reexpress~\eqref{eqn:98} as
$$
\left\{\begin{array}{l}
  A(-4k^2) =\pi/2K  \\
  B(-4k^2) =\frac{2E-K}{K}
\end{array}\right.
$$
This completes the proof of theorem~\ref{thm:main3}.\hfill$\Box$

\begin{rem}
\label{rem:Kr}
  The formula obtained in theorem~\ref{thm:main3} may be rewritten
  in the form given in~\cite[conjecture]{Kr} (hence providing a 
  proof of their conjecture), but the cohomology classes $K_i$
  appearing in~\cite[conjecture]{Kr} are the basic classes of
  $X$ and some multiples of them.
\end{rem}


\begin{thebibliography}{10}
\bibitem{basic}
V. Mu\~noz,
Basic classes for four-manifolds not of simple type,
{\tt math.DG/9811089}.
\end{thebibliography}
\end{document}